\documentclass[10pt,a4paper,leqno,english]{article}
\usepackage{graphicx}
\usepackage[T1]{fontenc}
\usepackage{babel}
\usepackage{enumerate,amsmath,amssymb,latexsym}

\newtheorem{example}{Example}[section]

\newtheorem{lemma}[example]{Lemma}
\newtheorem{proposition}[example]{Proposition}
\newtheorem{theorem}[example]{Theorem}

\newtheorem{corollary}[example]{Corollary}

\newenvironment{proof}[1][Proof]{\noindent\textit{#1.} }{\ \rule{0.5em}{0.5em}}
\newenvironment{keywords}[1][Keywords]{\noindent\textbf{#1.} }{\ \rule{0.5em}{0.5em}}
\newenvironment{M}[1][M.S.C. 2000]{\noindent\textbf{#1.} }{\ \rule{0.5em}{0.5em}}
\newenvironment{abs}[1][Abstract]{\noindent\textbf{#1.} }{\ \rule{0.5em}{0.5em}}

\begin{document}

\date{}
\title{Hypersurfaces of Lorentzian para-Sasakian Manifolds }
\author{Selcen Y\"{u}ksel Perkta\c{s}, Erol K\i l\i \c{c}, Sad\i k Kele\c{s}}
\maketitle

\abs{In this paper, we study the invariant and noninvariant
hypersurfaces of  $(1,1,1)$ almost contact manifolds, Lorentzian
almost paracontact manifolds and Lorentzian para-Sasakian manifolds,
respectively. We show that a noninvariant hypersurface of an
$(1,1,1)$ almost contact manifold admits an almost product
structure. We investigate hypersurfaces of affinely cosymplectic and
normal $(1,1,1)$ almost contact manifolds. It is proved that a
noninvariant hypersurface of a Lorentzian almost paracontact
manifold is an almost product metric manifold. Some necessary and
sufficient conditions have been given for a noninvariant
hypersurface of a Lorentzian para-Sasakian manifold to be locally
product manifold. We establish a Lorentzian para-Sasakian structure
for an invariant hypersurface of a Lorentzian para-Sasakian
manifold. Finally we give some examples for invariant and
noninvariant hypersurfaces of a
Lorentzian para-Sasakian manifold.\\
\M{53C25, 53C42, 53C50.}\\

\section{Introduction}

Hypersurfaces of an almost contact manifold have been studied by Blair \cite%
{Blair}, Eum \cite{Eum}, Goldberg and Yano \cite{Goldberg}, Ludden \cite%
{Ludden1} and the others.
 In 1970, Goldberg and Yano
\cite{Goldberg} defined noninvariant hypersurfaces of almost contact
manifolds. A hypersurface such
that the transform of a tangent vector of the hypersurface by the tensor $%
\varphi $ defining the almost contact structure is never tangent to
the
hypersurface is called a noninvariant hypersurface of almost contact manifold%
\cite{Goldberg}. The authors \cite{Goldberg} showed that a
noninvariant hypersurface of an almost contact manifolds admits an
almost complex structure and a distinguished 1-form induced by the
contact form of the manifold. They also investigated the
noninvariant hypersurface of an almost contact metric manifold.

In 1976, Sato \cite{Sato} studied a structure similiar to the almost
contact structure, namely almost paracontact structure. In
\cite{Adati}, T. Adati studied hypersurfaces of an almost
paracontact manifold. A. Bucki \cite{Bucki} considered hypersurfaces
of an almost r-paracontact Riemannian manifold. Some properties of
invariant hypersurfaces of an almost r-paracontact Riemannian
manifold were investigated in \cite{Bucki2} by A. Bucki and A.
Miernowski. Moreover, in \cite{Ion}, I. Mihai and K. Matsumoto
studied submanifolds of an almost r-paracontact Riemannian manifold
of P-Sasakian type.In \cite{Dube} the authors studied invariant and
noninvariant hypersurfaces of almost R-paracontact manifolds. Singh \cite{singh} defined $%
(e_{1},e_{2},r)$ almost contact structure as a generalization of
many known structures, which are obtained by taking particular values of $%
(e_{1},e_{2})$ and $r$ (see also \cite{mukutilk}). The study of
Lorentzian almost paracontact manifolds was initiated by Matsumoto
in 1989 \cite{F11}. Also he introduced the notion of Lorentzian
para-Sasakian ( for short, \textit{LP}-Sasakian ) manifold. I. Mihai
and R. Rosca \cite{H3} defined the
same notion independently and thereafter many authors \cite%
{H4,H5,mukut,mukut2} studied \textit{LP}-Sasakian manifolds and
their submanifolds.

In the present paper, we study invariant and noninvariant
hypersurfaces of $(1,1,1)$ almost contact manifolds, Lorentzian
almost paracontact manifolds and Lorentzian para-Sasakian manifolds,
respectively. We investigate the invariant hypersurfaces with two
different conditions: when the characteristic vector field $\xi$ is
everywhere tangent to the hypersurfaces and when the characteristic
vector field $\xi$ does not belong to the tangent hyperplane of the
hypersurfaces. Section 2 is devoted to preliminaries. In section 3
we show that a noninvariant hypersurface of a $(1,1,1)$ almost
contact manifold with the characteristic vector field $\xi$ nowhere
tangent to the hypersurface admits an almost product structure. In
section 4 we study hypersurfaces of affinely cosymplectic and normal
$(1,1,1)$ almost contact manifolds. In section 5 it is proved that a
noninvariant hypersurface of a Lorentzian almost paracontact
manifold is an almost product metric manifold. We also find a
necessary and sufficient condition for a noninvariant hypersurface
of a Lorentzian para-Sasakian manifold to be locally product
manifold. In section 5 we establish a Lorentzian para-Sasakian
structure for an invariant hypersurface of a Lorentzian
para-Sasakian manifold with the characteristic vector field $\xi$
tangent to the hypersurface. In the last section we give some
examples for invariant and noninvariant hypersurfaces of an
$(1,1,1)$ almost contact manifold, a Lorentzian almost paracontact
manifold and a Lorentzian para-Sasakian manifold.

\section{Preliminaries}

\setcounter{equation}{0}
\renewcommand{\theequation}{2.\arabic{equation}}

Let $\overline{M}$ be an n-dimensional differentiable manifold. If
there exist a tensor field $\varphi $ of type $(1,1)$, $r$-linearly
independent
vector fields $\xi _{\alpha }$ and $r$ 1-forms $\eta ^{\alpha }$ on $%
\overline{M}$ such that \cite{mukutilk}
\begin{eqnarray}
 \varphi (\xi
_{\alpha }) &=&0,\label{pre-1}\\
 {\varphi }^{2} &=&e_{1}I+e_{2}\eta
^{\alpha }\otimes \xi _{\alpha },\label{pre-2}
\end{eqnarray}
where $e_{1},$ $e_{2}$ take values $\pm 1$ independently, $I$
denotes the identity map of $\Gamma (T\overline{M})$ and $\otimes $
is the tensor product, then the structure $(\varphi,\xi _{\alpha
},\eta ^{\alpha })$ is said to be an almost
$(e_{1},e_{2})$-$r$-contact structure or in short $(e_{1},e_{2},r)$
ac structure and the manifold $\overline{M}$ with the
$(e_{1},e_{2},r)$ ac structure is called an $(e_{1},e_{2},r)$ ac
manifold.
\newpage

Let $\overline{M}$ be an $(e_{1},e_{2},r)$ ac manifold. Then the
following relations hold on $\overline{M}$ \cite{singh}:
\begin{eqnarray}
\eta ^{\alpha }\circ \varphi &=&0, \label{pre-3}\\
\eta ^{\alpha }(\xi _{\beta }) &=&-e_{1}e_{2}\delta _{\beta }^{\alpha }, \label{pre-4}\\
rank(\varphi ) &=&n-r\label{pre-5}.
\end{eqnarray}

Now, consider that $\overline{M}$ is a $(1,1,1)$ ac manifold. Then $\overline{M%
}$ admits a Lorentzian metric $\overline{g}$, such that%
\begin{eqnarray}
\eta (\overline{X}) &=&\overline{g}(\overline X,\xi ),\label{pre-6} \\
\overline{g}(\varphi \overline X,\varphi \overline Y)
&=&\overline{g}(\overline X,\overline Y)+\eta (\overline X)\eta
(\overline Y),\label{pre-7}
\end{eqnarray}%
for all $\overline X,\,\overline Y\in \Gamma (T\overline{M})$. In
this case $\overline{M}$ is said to admit a Lorentzian almost
paracontact structure $(\varphi ,\xi ,\eta ,\overline{g}).$ Then we
get
\begin{eqnarray}
\Phi (\overline X,\overline Y) &\equiv &\overline{g}(\overline
X,\varphi \overline Y)\equiv \overline{g}(\varphi
\overline X,\overline Y)\equiv \Phi (\overline Y,\overline X), \label{pre-8}\\
(\overline{\nabla }_{\overline X}\Phi )(\overline Y,\overline Z) &=&\overline{g}(\overline Y,(\overline{\nabla }%
_{\overline X}\varphi )Z)=(\overline{\nabla }_{\overline X}\Phi
)(\overline Z,\overline Y),\label{pre-9}
\end{eqnarray}%
where $\overline{\nabla }$ is the Levi-Civita connection with respect to $%
\overline{g}$. It is clear that Lorentzian metric $\overline{g}$
makes $\xi $ a timelike unit vector field, i.e, $\overline{g}(\xi
,\xi )=-1.$ The manifold $\overline{M}$ equipped with a Lorentzian
almost paracontact structure $(\phi ,\xi ,\eta ,\overline{g})$ is
called a Lorentzian almost \ paracontact manifold (for short
$LAP$-manifold) \cite{F11,F12}.

A Lorentzian almost \ paracontact manifold $\overline{M}$ endowed
with the structure $(\varphi ,\xi ,\eta ,\overline{g})$ is called a
Lorentzian
paracontact manifold ( for short \textit{LP}-manifold) \cite{F11} if%
\begin{eqnarray}
\Phi (\overline X,\overline Y)=\frac{1}{2}((\overline{\nabla }_{\overline X}\eta )\overline Y+(\overline{\nabla }%
_{\overline Y}\eta )\overline X).\label{pre-10}
\end{eqnarray}

A Lorentzian almost \ paracontact manifold $\overline{M}$ endowed
with the structure $(\varphi ,\xi ,\eta ,\overline{g})$ is called a
Lorentzian para
Sasakian manifold ( for short \textit{LP}-Sasakian) \cite{F11} if%
\begin{eqnarray}
(\overline{\nabla }_{\overline X}\varphi )\overline Y=\eta
(\overline Y)\overline X+\overline{g}(\overline X,\overline Y)\xi
+2\eta (\overline X)\eta (\overline Y)\xi .\label{pre-11}
\end{eqnarray}%
We note that in a \textit{LP}-Sasakian manifold the $1$-form $\eta $
is closed.

Let $\overline{M}\times R$ be a product manifold, where
$\overline{M}$ is an
$(1,1,1)$ ac manifold. The tensor field $J^{\prime }$ of type $(1,1)$ on $%
\overline{M}\times R$ defined by
\begin{eqnarray}
J^{\prime }(\overline X,f\frac{d}{dt})=(\varphi \overline X-f\xi
,\eta (\overline X)\frac{d}{dt}),\label{pre-12}
\end{eqnarray}%
where $f$ is a $C^{\infty }$ real-valued function \ and $\overline X\in \Gamma (T%
\overline{M}),$ satisfies $J^{\prime 2}=I$ and thus provides an
almost product structure on $\overline{M}\times R.$ If the induced
almost product structure on $\overline{M}\times R$ is integrable
then the $(1,1,1)$ ac structure on \ $\overline{M}$ is said to be
normal \cite{mukutilk}. Since the vanishing of the Nijenhius tensor
$[J^{\prime },J^{\prime }]$ is a necessary and sufficient condition
for integrability, the condition of the normality in terms of
Nijenhius tensor $[\varphi ,\varphi ]$ of $\varphi $ is (see
\cite{mukutilk})
\begin{eqnarray}
\lbrack \varphi ,\varphi ]+d\eta \otimes \xi =0,\label{pre-13}
\end{eqnarray}%
where
\begin{equation}
\lbrack \varphi ,\varphi ](\overline X,\overline Y)=[\varphi
\overline X,\varphi \overline Y]-\varphi \lbrack \varphi \overline
X,\overline Y]-\varphi \lbrack \overline X,\varphi \overline
Y]+\varphi ^{2}[\overline X,\overline Y],\label{pre-14}
\end{equation}
for all $\overline X,\overline Y\in \Gamma(T\overline M)$.
\section{Non-Invariant Hypersurfaces of $(1,1,1)$ ac\\
 Manifolds}

\setcounter{equation}{0}
\renewcommand{\theequation}{3.\arabic{equation}}

Let $\overline{M}$ is a $(1,1,1)$ ac manifold. Consider an
$(n-1)$-dimensional
manifold $M$ imbedded in $\overline{M}$ with the imbedding map%
\begin{equation*}
i:M\rightarrow \overline{M},
\end{equation*}%
and assume that for each $p\in M$ the vector $\xi _{i(p)}$ does not
belong to the tangent hyperplane of the hypersurface. Then we have
\begin{equation}
\varphi i_{\ast }X=i_{\ast }JX+\alpha (X)\xi ,\label{noninvariant-1}
\end{equation}%
where $J$ is a tensor field of type $(1,1)$ , $\alpha $ is a
$1$-form on $M$
and \ $i_{\ast }$ is the differential of the immersion $i$ of $M$ into $%
\overline{M}.$ If $\alpha \neq 0$, then the submanifold $i(M)$ is
called a
noninvariant hypersurface of $\overline{M}$. On the other hand, if the $1$%
-form $\alpha $ vanishes, $i(M)$ is called an invariant hypersurface of $%
\overline{M}$ (see \cite{Goldberg}). A hypersurface may, of course,
be neither invariant nor
noninvariant. Throughout this section, unless specified otherwise $%
i(M)$ will be a noninvariant hypersurface of the $(1,1,1)$ ac manifold $%
\overline{M}.$

\begin{theorem}
If $M$ is a noninvariant hypersurface of a $(1,1,1)$ ac manifold
$\overline{M}$ with $\xi$ nowhere tangent to $M$, then $M$ admits an
almost product structure.
\end{theorem}

\begin{proof}
By applying $\varphi $ to (\ref{noninvariant-1}) and using (\ref{pre-1})-(\ref{pre-4}), we have%
\begin{equation}
i_{\ast }X+\eta (i_{\ast }X)\xi =i_{\ast }(J^{2}X)+\alpha (JX)\xi
.\label {noninvariant-2}
\end{equation}%
Then from (\ref{noninvariant-1}), we get%
\begin{equation*}
J^{2}X=X
\end{equation*}%
and%
\begin{equation}
\alpha (JX)=\eta (i_{\ast }X)=i^{\ast }(\eta
X),\label{noninvariant-3}
\end{equation}%
where $X\in \Gamma (TM)$ and $i^{\ast }$ is the dual map of $i_{\ast
}.$ Thus $J$ acts as an almost product structure on $M$.
\end{proof}

If we define a $1$-form $C\alpha $ on $M$ by $C\alpha (X)=\alpha
(JX)$ then from (\ref{noninvariant-3}) we can write
\begin{equation*}
C\alpha =i^{\ast }\eta .
\end{equation*}%
Thus, the hypersurface $M$ admits a $1$-form $\alpha $ whose
vanishing means
that the tangent hyperplane of the hypersurface is invariant under $\varphi $%
.

Now, let $\overline{\nabla }$ be a symmetric affine connection on $\overline{%
M}$ and define an affine connection $\nabla $ on $M$ with respect to
the affine normal $\xi $ by
\begin{equation}
\overline{\nabla }_{i_{\ast }X}i_{\ast }Y=i_{\ast }\nabla
_{X}Y+h(X,Y)\xi ,\label{noninvariant-4}
\end{equation}%
where $h$ is a symmetric tensor field of type (0,2) on $M$ which is
called the second fundamental form of $M$ with respect to $\xi$.

Suppose that the $(1,1,1)$ ac structure is normal. Then, the torsion
field $S$
of type $(1,2)$ on $M$ which is defined by%
\begin{equation}
S(\overline{X},\overline{Y})=[\varphi \overline{X},\varphi \overline{Y}%
]-\varphi \lbrack \varphi \overline{X},\overline{Y}]-\varphi \lbrack
\overline{X},\varphi \overline{Y}]+\varphi ^{2}[\overline{X},\overline{Y}%
]+d\eta (\overline{X},\overline{Y})\xi ,\label{noninvariant-5}
\end{equation}%
for all $\overline{X},\overline{Y}\in \Gamma (T\overline{M}),$
vanishes. By
taking $\overline{Y}=\xi $ in (\ref{noninvariant-5}), we get%
\begin{equation*}
L_{\xi }\varphi =0\text{ \ \ \ and \ \ \ }L_{\xi }\eta =0,
\end{equation*}%
where $L_{\xi }$ is the Lie derivative operator with respect to $\xi
.$ From
(\ref{noninvariant-5}) the tensor field $S$ is also expressed by%
\begin{eqnarray}
S(\overline{X},\overline{Y}) &=&\overline{\nabla }_{\varphi \overline{X}%
}(\varphi \overline{Y})-\overline{\nabla }_{\varphi
\overline{Y}}(\varphi
\overline{X})-\varphi (\overline{\nabla }_{\varphi \overline{X}}\overline{Y}-%
\overline{\nabla }_{\overline{Y}}(\varphi \overline{X})) \nonumber\\
&&-\varphi (\overline{\nabla }_{\overline{X}}(\varphi \overline{Y})-%
\overline{\nabla }_{\varphi \overline{Y}}\overline{X})+\varphi ^{2}(%
\overline{\nabla }_{\overline{X}}\overline{Y}-\overline{\nabla }_{\overline{Y%
}}\overline{X}) \label{noninvariant-6}\\
&&+(\overline{\nabla }_{\overline{X}}\eta (\overline{Y})-\overline{\nabla }_{%
\overline{Y}}\eta (\overline{X})-\eta
([\overline{X},\overline{Y}]))\xi \nonumber,
\end{eqnarray}
or
\begin{eqnarray}
S(\overline{X},\overline{Y}) &=&(\overline{\nabla }_{\varphi \overline{X}%
}\varphi )\overline{Y}-(\overline{\nabla }_{\varphi \overline{Y}}\varphi )%
\overline{X}+\varphi (\overline{\nabla }_{\overline{Y}}\varphi )\overline{X}%
-\varphi (\overline{\nabla }_{\overline{X}}\varphi )\overline{Y} \label{noninvariant-7}\\
&&+[(\overline{\nabla }_{\overline{X}}\eta )\overline{Y}-(\overline{\nabla }%
_{\overline{Y}}\eta )\overline{X}]\xi .  \notag
\end{eqnarray}
By using (\ref{noninvariant-1}) and (\ref{noninvariant-4}), we
obtain
\begin{eqnarray}
S(i_{\ast }X,i_{\ast }Y) &=&i_{\ast }[J,J](X,Y)+L_{\xi }\varphi
\{\alpha
(X)i_{\ast }Y-\alpha (Y)i_{\ast }X\} \\
&&+\{d\alpha (JX,Y)+d\alpha (X,JY)-2i^{\ast }\eta ([X,Y])\}\xi
\notag
\end{eqnarray}%
Therefore, we have

\begin{theorem}
A noninvariant hypersurface of a normal $(1,1,1)$ ac manifold
$\overline{M}$ is a locally product manifold which has a $1$-form
$\alpha =C^{-1}i^{\ast
}\eta $ such that its differential satisfies%
\begin{equation}
d\alpha (JX,Y)+d\alpha (X,JY)=2C\alpha
([X,Y]).\label{noninvariant-9}
\end{equation}
\end{theorem}

\begin{corollary}
An invariant hypersurface of a $(1,1,1)$ ac manifold is an almost
product manifold. If the $(1,1,1)$ ac manifold is normal, then the
almost product structure is integrable.
\end{corollary}

\begin{theorem}
Let $\xi$ be an infinitesimal automorphism of the $(1,1,1)$ ac manifold $%
\overline{M}$. If, for every noninvariant hypersurface, the induced
almost
product structure $J$ is integrable and the differential of the induced $1$%
-form $\alpha =C^{-1}i^{\ast }\eta $ satisfies
(\ref{noninvariant-9}) then $\overline{M}$ is normal.
\end{theorem}

\section{Hypersurfaces of affinely cosymplectic and\\
 normal $(1,1,1)$ ac
manifolds}

\setcounter{equation}{0}
\renewcommand{\theequation}{4.\arabic{equation}}

Let $\overline{M}$ be a $(1,1,1)$ ac manifold with a symmetric
affine connection $\overline{\nabla }$ and $\nabla $ denotes the
induced connection
on the noninvariant hypersurface $M$. If we write%
\begin{equation}
(\nabla _{X}i_{\ast })Y=\overline{\nabla }_{i_{\ast }X}i_{\ast
}Y-i_{\ast }(\nabla _{X}Y),\label{affine-1}
\end{equation}%
then the Gauss and Weingarten equations are%
\begin{equation}
(\nabla _{X}i_{\ast })Y=h(X,Y)\xi ,\text{ \ \ \
}h(X,Y)=h(Y,X),\label{affine-2}
\end{equation}%
and
\begin{equation}
\overline{\nabla }_{i_{\ast }X}\xi =-i_{\ast }AX+w(X)\xi
,\label{affine-3}
\end{equation}%
where $h$ and $A$ are the second fundamental tensors of type $(0,2)$ and $%
(1,1)$, respectively of $M$ with respect to $\xi ,$ and $w$ is a 1-form on $%
M $ defining the connection on the affine normal bundle.

By using (\ref{noninvariant-1}), (\ref{affine-2}) and (\ref{affine-3}) we get%
\begin{eqnarray}
(\overline{\nabla }_{i_{\ast }X}\varphi )i_{\ast }Y &=&\overline{\nabla }%
_{i_{\ast }X}\varphi i_{\ast }Y-\varphi (\overline{\nabla }_{i_{\ast
}X}i_{\ast }Y)  \notag \\
&=&[h(X,JY)+(\nabla _{X}\alpha )(Y)+w(X)\alpha (Y)]\xi \label{affine-4}\\
&&+i_{\ast }[(\nabla _{X}J)Y-\alpha (Y)AX] .  \notag
\end{eqnarray}%
Then we will investigate the following two cases:\\

\textbf{Case I : } Let $\overline{M}$ be an affinely cosymplectic
$(1,1,1)$ ac manifold, that is, $\overline{M}$ be a $(1,1,1)$ ac
manifold with a symmetric
affine connection $\overline{\nabla }$ such that%
\begin{equation}
\overline{\nabla }\varphi =0\text{, \ \ \ \ \ \ \ }\overline{\nabla
}\eta =0.\label{affine-5}
\end{equation}%
From (\ref{noninvariant-7}) we can easily see that an affinely
cosymplectic $(1,1,1)$ ac manifold is normal. Also by using
(\ref{pre-1}) and (\ref{pre-2}), we can show that (\ref{affine-5})
implies that
\begin{equation*}
\overline{\nabla }\xi =0.
\end{equation*}%
Therefore, by (\ref{affine-3}), we have
\begin{equation*}
AX=0\text{ \ \ and \ }w(X)=0.
\end{equation*}%
Moreover, since $\overline{\nabla }\varphi =0$ then from
(\ref{affine-4}) we have
\begin{equation*}
\nabla J=0
\end{equation*}%
and%
\begin{equation*}
(\nabla _{X}\alpha )(Y)=-h(X,JY).
\end{equation*}

\textbf{Case II : } Let $\overline{M}$ be a normal $(1,1,1)$ ac
manifold such that $\varphi =\overline{\nabla }\xi .$ Then by using
(\ref{noninvariant-1}) and (\ref{affine-3}), we
have%
\begin{equation*}
i_{\ast }JX+\alpha (X)\xi =-i_{\ast }AX+w(X)\xi
\end{equation*}%
that \ is%
\begin{equation*}
J=-A
\end{equation*}%
and%
\begin{equation*}
\alpha =w.
\end{equation*}

If $AX=0,$ for all $X\in \Gamma (TM)$, then from (\ref{affine-3}) it is obvious that $%
\overline{\nabla }_{i_{\ast }X}\xi $ and $\xi $ are propotional. So
affine
normals are parallel along the hypersurface. In this case, the hypersurface $%
M$ is said to be totally flat.

\begin{proposition}
Let $M$ be a noninvariant hypersurface of an affinely cosymplectic
$(1,1,1)$
ac manifold. Then $M$ is totally flat and%
\begin{eqnarray*}
\nabla J &=&0, \\
(\nabla _{X}\alpha )(Y) &=&-h(X,JY), \\
w &=&0.
\end{eqnarray*}
\end{proposition}

\begin{corollary}
Let $M$ be an invariant hypersurface of an affinely cosymplectic
$(1,1,1)$ ac manifold. Then
\begin{eqnarray*}
\nabla J &=&0, \\
h &=&0, \\
w &=&0.
\end{eqnarray*}
\end{corollary}

\begin{proposition}
Let $M$ be a noninvariant hypersurface of a normal $(1,1,1)$ ac
manifold such
that $\varphi =\overline{\nabla }\xi .$ Then%
\begin{equation*}
J=-A
\end{equation*}%
and%
\begin{equation*}
\alpha =w.
\end{equation*}
\end{proposition}

\section{Hypersurfaces of Lorentzian almost\\
 paracontact manifolds}

\setcounter{equation}{0}
\renewcommand{\theequation}{5.\arabic{equation}}

A $(1,1,1)$ ac manifold $\overline{M}$ admitting a Lorentzian metric $%
\overline{g}$ such that%
\begin{eqnarray}
\overline{g}(\overline{X},\xi ) &=&\eta (\overline{X})\label{lap-1})\\
\overline{g}(\overline{X},\varphi \overline{Y}) &\equiv
&\overline{g}(\varphi \overline{X},\overline{Y}),\label{lap-2}
\end{eqnarray}
where $\overline{X},\overline{Y}\in\Gamma(T\overline{M})$, is called Lorentzian almost paracontact manifold and denoted by $(\overline{M%
},\varphi ,\eta ,\overline{g}).$

\begin{proposition}
Let $(M,J,\alpha ,g)$ be a noninvariant hypersurface of $(\overline{M}%
,\varphi ,\eta ,\overline{g})$ where $g$ is the induced metric on
$M$, that is, $i^{\ast }\overline{g}=g.$ Then the hypersurface
$(M,J,\alpha ,g)$ admits an almost product metric
\begin{eqnarray}
G=g+\alpha\otimes \alpha.\label{lap-3}
\end{eqnarray}
\end{proposition}

\begin{proof}
From (\ref{lap-2}), we can write
\begin{equation}
\overline{g}(\varphi i_{\ast }X,i_{\ast }Y)=\overline{g}(\varphi
i_{\ast }X,i_{\ast }Y).\label{lap-4}
\end{equation}%
By using (\ref{pre-1}) in (\ref{lap-4}), we obtain%
\begin{equation}
\overline{g}(i_{\ast }JX,i_{\ast }Y)+\alpha (X)\eta (i_{\ast }Y)=\overline{g}%
(i_{\ast }X,i_{\ast }JY)+\alpha (Y)\eta (i_{\ast }X).\label{lap-5}
\end{equation}%
The induced metric $g$ on $(M,J,\alpha )$ can be defined by
\begin{equation*}
g(X,Y)=\overline{g}(i_{\ast }X,i_{\ast }Y).
\end{equation*}%
So if we use (\ref{noninvariant-3}) and (\ref{lap-4}) in (\ref{lap-5}), then we have%
\begin{equation*}
g(JX,Y)+\alpha (X)C\alpha (Y)=g(X,JY)+\alpha (Y)C\alpha (X),
\end{equation*}%
that is,
\begin{equation*}
(g+\alpha \otimes \alpha )(JX,Y)=(g+\alpha \otimes \alpha )(X,JY).
\end{equation*}
If we denote $g+\alpha \otimes \alpha$ by $G$, the proof is
completed.
\end{proof}

\begin{corollary}
A noninvariant hypersurface of a Lorentzian almost paracontact
manifold is an almost product metric manifold.
\end{corollary}

Now, let define 2-forms
\begin{equation*}
\Phi (\overline{X},\overline{Y})=\overline{g}(\varphi \overline{X},\overline{%
Y})\text{, \ \ }\overline{X},\overline{Y}\in \Gamma (T\overline{M})
\end{equation*}%
and%
\begin{equation*}
\Omega (X,Y)=G(JX,Y),\text{ \ \ \ \ \ }X,Y\in \Gamma (TM).
\end{equation*}%
$\Phi $ $\ $and $\ \Omega $ are called the fundamental forms of the
Lorentzian almost paracontact manifold $(\overline{M},\varphi ,\eta ,%
\overline{g})$ and the submanifold $(M,J,G)$ of $\overline M$,
respectively. Then we have

\begin{lemma}
Let $\Phi $ $\ $and $\ \Omega $ are the fundamental forms of \ $(\overline{M}%
,\varphi ,\eta ,\overline{g})$ and $(M,J,\alpha,G)$, respectively. Then%
\begin{equation}
i^{\ast }\Phi =\Omega -C\alpha \wedge \alpha .\label{lap-6}
\end{equation}
\end{lemma}

\begin{proof}
For $X,Y\in \Gamma (TM)$, by using definitions of the fundamental
forms, (\ref{noninvariant-1}) and (\ref{lap-3}), we get
\begin{eqnarray*}
\Phi (i_{\ast }X,i_{\ast }Y)=\Omega (X,Y)-(C\alpha \wedge \alpha
)(X,Y).
\end{eqnarray*}
Hence, we obtain
\begin{equation*}
i^{\ast }\Phi (X,Y)=(\Omega -C\alpha \wedge \alpha)(X,Y).
\end{equation*}
\end{proof}

\begin{theorem}
Let $(M,J,\alpha ,G)$ be a noninvariant hypersurface of the
Lorentzian para-Sasakian manifold $(\overline{M},\varphi ,\eta
,\overline{g})$. Then
\end{theorem}

\begin{description}
\item[(a)] \ $J=-A$,

\item[(b)] $\alpha=w$.
\end{description}

\begin{proof}
Since $(\overline{M},\varphi ,\eta ,\overline{g})$ is a Lorentzian
para-Sasakian manifold, we have%
\begin{equation*}
\overline{\nabla }_{i_{\ast }X}\xi =\varphi i_{\ast }X.
\end{equation*}%
By using (\ref{affine-3}) and (\ref{noninvariant-1}), we get
\begin{equation*}
-i_{\ast }AX+w(X)\xi =i_{\ast }JX+\alpha (X)\xi
\end{equation*}%
which gives%
\begin{equation*}
J=-A\text{ \ \ and \ }\alpha =w.
\end{equation*}
\end{proof}

\begin{theorem}
If $M$ is a noninvariant hypersurface of a Lorentzian para-Sasakian
manifold $(\overline{M},\varphi ,\eta ,\overline{g})$, then

\begin{description}
\item[(a)] $(\nabla _{X}J)(Y)=\alpha (Y)JX-C\alpha (Y)X,$

\item[(b)] $\overline{g}(i_{\ast }X,i_{\ast }Y)+2C\alpha (X)C\alpha
(Y)=h(X,JY)+(\nabla _{X}\alpha )(Y)+\alpha (X)\alpha (Y).$
\end{description}
\end{theorem}

\begin{proof}
By using (\ref{noninvariant-1}) and (\ref{affine-1}) we obatin%
\begin{eqnarray}
(\overline{\nabla }_{i_{\ast }X}\varphi )(i_{\ast }Y) &=&[i_{\ast
}(\nabla
_{X}J)(Y)+\alpha (Y)i_{\ast }JX]\label{lap-7} \\
&&+[h(X,JY)+(\nabla _{X}\alpha )(Y)+\alpha (X)\alpha (Y)]\xi .
\notag
\end{eqnarray}%
On the other hand, since $(\overline{M},\varphi ,\eta
,\overline{g})$ is a
Lorentzian para-Sasakian manifold, from (\ref{pre-11}) we also have%
\begin{equation}
(\overline{\nabla }_{i_{\ast }X}\varphi )(i_{\ast }Y)=\eta (i_{\ast
}Y)i_{\ast }X+\overline{g}(i_{\ast }X,i_{\ast }Y)\xi +2\eta (i_{\ast
}X)\eta (i_{\ast }Y)\xi .\label{lap-8}
\end{equation}%
By considering $C\alpha =i^{\ast }\eta $ and equating the components
of (\ref{lap-7}) and (\ref{lap-8}), we get (a) and (b) in the
theorem.
\end{proof}

As an immediate consequence we have the following:

\begin{corollary}
Let $M$ be a noninvariant hypersurface of the Lorentzian
para-Sasakian manifold $(\overline{M},\varphi ,\eta ,\overline{g})$
with the induced almost product structure $J.$ Then $M$ is a locally
product manifold if and only if
\begin{equation}
\alpha (Y)JX=\alpha (JY)X.\label{lap-9}
\end{equation}
\end{corollary}

Now, let $\overline{M}$ be a (1,1,1) ac manifold and $M$ be an
invariant hypersurface of $\overline{M}.$ Assume that for each $p\in
M$ the vector $\xi_{i(p)}$ belongs to the tangent hyperplane of the
hypersurface. For an invariant hypersurface of a $(1,1,1)$ ac
manifold we can write
\begin{equation}
\varphi i_{\ast }X=i_{\ast }\psi X,\label{lap-10}
\end{equation}%
where $\psi $ is a tensor of type (1,1) on the hypersurface $M$ and
$X\in \Gamma (TM).$ Applying $\varphi $ to the both sides of the
equation (\ref{lap-10}),
we get%
\begin{equation}
i_{\ast }\psi ^{2}X=\varphi ^{2}i_{\ast }X=i_{\ast }X+\eta (i_{\ast
}X)\xi .\label{lap-11}
\end{equation}%
If we denote
\begin{equation}
i_{\ast }\xi ^{\ast }=\xi\label{lap-12}
\end{equation}%
and
\begin{equation}
\eta ^{\ast }(X)=\eta (i_{\ast }X),\label{lap-13}
\end{equation}%
then we have%
\begin{equation}
\psi ^{2}X=X+\eta ^{\ast }(X)\xi ^{\ast }.\label{lap-14}
\end{equation}%
Furthermore,%
\begin{eqnarray}
\eta ^{\ast }(\psi X) &=&\eta (i_{\ast }\psi X)=\eta (\varphi
i_{\ast }X)=0,\label{lap-15}
\\
\eta ^{\ast }(\xi ^{\ast }) &=&\eta (i_{\ast }\xi ^{\ast })=\eta
(\xi )=-1\label{lap-16}
\end{eqnarray}%
and%
\begin{equation*}
i_{\ast }\psi \xi ^{\ast }=\varphi i_{\ast }\xi ^{\ast }=\varphi \xi
=0,
\end{equation*}%
that is
\begin{equation}
\psi \xi ^{\ast }=0.\label{lap-17}
\end{equation}%
Thus we have,

\begin{theorem}
Let $M$ be an invariant hypersurface of a $(1,1,1)$ ac manifold $(\overline{M}%
,\varphi ,\eta ,\xi )$ and $\xi\in\Gamma(TM)$. Then $M$ is a
$(1,1,1)$ ac manifold with the structure
$(\psi ,\xi ^{\ast },\eta ^{\ast })$ where $i_{\ast }\xi ^{\ast }=\xi$ and $%
\eta ^{\ast }(X)=\eta (i_{\ast }X)$, for all $X\in \Gamma (TM)$.
\end{theorem}

\begin{theorem}
Let $M$ be an invariant hypersurface of a $(1,1,1)$ ac manifold $(\overline{M}%
,\varphi ,\eta ,\xi )$ with $\xi\in\Gamma(TM)$. If $\overline{M}$ is
normal, then $M$ is also normal.
\end{theorem}

\begin{proof}
By using (\ref{noninvariant-5}), we can write%
\begin{eqnarray}
S(i_{\ast }X,i_{\ast }Y) &=&[\varphi ,\varphi ](i_{\ast }X,i_{\ast
}Y)+d\eta
(i_{\ast }X,i_{\ast }Y)\xi  \notag \\
&=&[\varphi i_{\ast }X,\varphi i_{\ast }Y]-\varphi \lbrack \varphi
i_{\ast
}X,i_{\ast }Y]-\varphi \lbrack i_{\ast }X,\varphi i_{\ast }Y]\label{lap-18} \\
&&+\varphi ^{2}[i_{\ast }X,i_{\ast }Y]+d\eta (i_{\ast }X,i_{\ast
}Y)\xi . \notag
\end{eqnarray}%
for all $X,Y\in \Gamma(TM)$. If we use (\ref{lap-10}),
(\ref{lap-12}) and (\ref{lap-13}) in (\ref{lap-18}),
we get%
\begin{eqnarray*}
S(i_{\ast }X,i_{\ast }Y) &=&i_{\ast }\psi ^{2}[X,Y]+[i_{\ast }\psi
X,i_{\ast
}\psi Y]-i_{\ast }\psi \lbrack X,\psi Y]-i_{\ast }\psi \lbrack \psi X,Y] \\
&&+\{(i_{\ast }X)(\eta ^{\ast }(Y))-(i_{\ast }Y)(\eta ^{\ast
}(X))-\eta
^{\ast }([X,Y])\}i_{\ast }\xi ^{\ast } \\
&=&i_{\ast }\{[\psi ,\psi ](X,Y)+d\eta ^{\ast }(X,Y)\xi ^{\ast }\}.
\end{eqnarray*}%
Hence, we have the assertion of the theorem.
\end{proof}

\begin{theorem}
Let $M$ be an invariant hypersurface of a Lorentzian almost
paracontact manifold $(\overline{M},\varphi ,\eta ,\overline{g})$
where $\xi\in\Gamma(TM)$. Then $M$ is also a Lorentzian almost
paracontact manifold.
\end{theorem}

\begin{proof}
From Theorem 5.7. it follows that an invariant hypersurface $M$ of $%
\overline{M}$ is a $(1,1,1)$ ac manifold with the structure $(\psi
,\xi ^{\ast
},\eta ^{\ast }).$ Let $g^{\ast }$ be the induced metric on $M$. Then we have%
\begin{equation}
g^{\ast }(\psi X,\psi Y)=\overline{g}(i_{\ast }\psi X,i_{\ast }\psi Y)=%
\overline{g}(\varphi i_{\ast }X,\varphi i_{\ast }Y).\label{lap-19}
\end{equation}%
Since $\overline{M}$ is a Lorentzian almost paracontact manifold,
then by
using (\ref{lap-13}) in (\ref{lap-19}) we get%
\begin{equation}
g^{\ast }(\psi X,\psi Y)=g^{\ast }(X,Y)+\eta ^{\ast }(X)\eta ^{\ast
}(Y).\label{lap-20}
\end{equation}%
Moreover,%
\begin{equation}
g^{\ast }(X,\xi ^{\ast })=\overline{g}(i_{\ast }X,i_{\ast }\xi
^{\ast })=\eta (i_{\ast }X)=\eta ^{\ast }(X),\label{lap-21}
\end{equation}%
which completes the proof.
\end{proof}

\begin{theorem}
Let $(\overline{M},\varphi ,\eta ,\overline{g})$ be a Lorentzian
para Sasakian manifold$.$Then an invariant hypersurface with
$\xi\in\Gamma(TM)$ of $\overline{M}$ is also a Lorentzian para
Sasakian manifold.
\end{theorem}

\begin{proof}
Let $\overline{M}$ be a Lorentzian para-Sasakian manifold. Then we
have
\begin{equation*}
\overline{\nabla }_{i_{\ast }X}\xi =\varphi i_{\ast }X,
\end{equation*}
where $\overline{\nabla }$ is a Levi-Civita connection with respect to $%
\overline{g}.$ From (\ref{lap-10}) and (\ref{lap-12}), we can write%
\begin{equation*}
\overline{\nabla }_{i_{\ast }X}i_{\ast }\xi ^{\ast }=i_{\ast }\psi
X.
\end{equation*}%
By using (\ref{noninvariant-4}) in the last equation, we obtain%
\begin{equation*}
i_{\ast }\nabla _{X}\xi ^{\ast }+h(X,\xi ^{\ast })N=i_{\ast }\psi X,
\end{equation*}%
where $\nabla $ is the induced connection on $M$ and $N$ is normal
to $M.$
If we consider normal and tangent components of above equation we get%
\begin{eqnarray*}
\nabla _{X}\xi ^{\ast } &=&\psi X, \\
h(X,\xi ^{\ast }) &=&0.
\end{eqnarray*}
Since $\overline{M}$ be a Lorentzian para Sasakian manifold from
(\ref{pre-11}), we
have%
\begin{equation}
(\overline{\nabla }_{i_{\ast }X}\varphi )i_{\ast }Y=\eta (i_{\ast
}Y)i_{\ast }X+\overline{g}(i_{\ast }X,i_{\ast }Y)\xi +2\eta (i_{\ast
}X)\eta (i_{\ast }Y)\xi ,\label{lap-22}
\end{equation}%
for all $X,Y\in \Gamma (TM).$ By using (\ref{lap-10}),
(\ref{lap-12}) and (\ref{lap-13}) in (\ref{lap-22}),
we obtain%
\begin{equation}
(\overline{\nabla }_{i_{\ast }X}\varphi )i_{\ast }Y=i_{\ast }\{\eta
^{\ast }(Y)X+\overline{g}(X,Y)\xi ^{\ast }+2\eta ^{\ast }(X)\eta
^{\ast }(Y)\xi ^{\ast }\}.\label{lap-23}
\end{equation}%
On the other hand, from (\ref{noninvariant-4}) and (\ref{lap-10}),
one can get
\begin{eqnarray}
(\overline{\nabla }_{i_{\ast }X}\varphi )i_{\ast }Y &=&\overline{\nabla }%
_{i_{\ast }X}\varphi i_{\ast }Y-\varphi (\overline{\nabla }_{i_{\ast
}X}i_{\ast }Y)  \notag \\
&=&\overline{\nabla }_{i_{\ast }X}i_{\ast }\psi Y-\varphi (i_{\ast
}\nabla
_{X}Y+h(X,Y)N)  \notag \\
&=&i_{\ast }(\nabla _{X}\psi Y-\psi (\nabla _{X}Y))+h(X,\psi
Y)N-h(X,Y)\varphi N,\label{lap-24}
\end{eqnarray}%
where $\nabla $ is the induced connection on $M$ and $N$ is normal
to $M.$By
equating right hand sides of equations (\ref{lap-23}) and (\ref{lap-24}), we have%
\begin{equation*}
(\nabla _{X}\psi )Y=\eta ^{\ast }(Y)X+\overline{g}(X,Y)\xi ^{\ast
}+2\eta ^{\ast }(X)\eta ^{\ast }(Y)\xi ^{\ast }.
\end{equation*}%
This completes the proof.
\end{proof}

\section{Examples}
 \setcounter{equation}{0}
\renewcommand{\theequation}{6.\arabic{equation}}
\begin{example}
Let $ \overline{M}$, be the 5-dimensional real number space with a
coordinate system $(x,y,z,t,s)$. Defining
\[
\eta =ds-dx-dz\ ,\qquad \xi =-\,\frac{\partial }{\partial s}\ ,
\]%
\[
\varphi \left( \frac{\partial }{\partial x}\right) =-\,\frac{\partial }{%
\partial x}-\frac{\partial }{\partial s}\ ,\qquad \varphi \left( \frac{%
\partial }{\partial y}\right) =-\,\frac{\partial }{\partial y}\ ,
\]%
\[
\varphi \left( \frac{\partial }{\partial z}\right) =-\,\frac{\partial }{%
\partial z}-\frac{\partial }{\partial s}\ ,\qquad \varphi \left( \frac{%
\partial }{\partial t}\right) =-\,\frac{\partial }{\partial t}\ ,\qquad
\varphi \left( \frac{\partial }{\partial s}\right) =0\ ,
\]%
the set $(\varphi ,\xi ,\eta )$\ becomes a $(1,1,1)$ ac structure in
$ \overline{M}$.

Let $M_{1}$ be a hypersurface of $ \overline{M}$ which is given by
$s=x$ with the imbedding map $i:M_{1}\rightarrow \overline{M} $ .
Then
\[
\{u_{1}=(1,0,0,0,1),u_{2}=(0,1,0,0,0),u_{3}=(0,0,1,0,0),u_{4}=(0,0,0,1,0)\}
\]
is a local basis for the tangent hyperplane of $M_{1}$  and
$N_{1}=(1,0,0,0,-1)$ is the normal vector field of the hypersurface.
It is obvious that the characteristic vector field $\xi_{i(p)}$,
$p\in M_{1}$, does not belong to the tangent hyperplane of $M_{1}$.
A tangent vector field of the hypersurface can be written by $X
\equiv i_{*}X=f_{1}u_{1}+f_{2}u_{2}+f_{3}u_{3}+f_{4}u_{4}$ for some
smooth functions $f_{i}$, $1\leq i\leq 4$, on $M$. Then we have
\begin{eqnarray*}
\varphi i_{*}X=-f_{1}u_{1}-f_{2}u_{2}-f_{3}u_{3}-f_{4}u_{4}+f_{3}\xi
\end{eqnarray*}
which shows that $M_{1}$ is a noninvariant hypersurface of
$\overline M$.

Now let us consider the hypersurface $M_{2}$ of the $(1,1,1)$ ac
manifold $ \overline{M}$ defining by $x=y$ and let
$i:M_{2}\rightarrow \overline{M} $ be the imbedding map of $M_{2}$
into $M$. In this case the set
\[
\{v_{1}=(1,1,0,0,0),v_{2}=(0,0,1,0,0),v_{3}=(0,0,0,1,0),v_{4}=(0,0,0,0,1)\}
\]
is a local basis for the tangent hyperplane and $N_{2}=(1,-1,0,0,0)$
is the normal vector field of $M_{2}$. The characteristic vector
field belongs to the tangent hyperplane of the hypersurface. For any
tangent vector field $X \equiv
i_{*}X=h_{1}v_{1}+h_{2}v_{2}+h_{3}v_{3}+h_{4}v_{4}$ of the
hypersurface we have
\begin{eqnarray}
\varphi i_{*}X=-h_{1}v_{1}-h_{2}v_{2}-f_{3}v_{3}+(h_{1}+h_{2})\xi
\label {ex-1}
\end{eqnarray}
where $h_{i}$, $1\leq i\leq 4$, are some smooth functions on
$M_{2}$. From (\ref{ex-1}) we see that that $M_{2}$ is an invariant
hypersurface of $\overline M$.
\end{example}

\begin{example}
Let $ \overline{M}$ be the 5-dimensional real number space with a
coordinate system $(x,y,z,t,s)$. In $ \overline{M}$ we define
\[
\eta =ds-dx\ ,\qquad \xi =-\,\frac{\partial }{\partial s}\ ,
\]%
\[
\varphi \left( \frac{\partial }{\partial x}\right) =\frac{\partial }{%
\partial x}+\frac{\partial }{\partial s}\ ,\qquad \varphi \left( \frac{%
\partial }{\partial y}\right) =\frac{\partial }{\partial y}\ ,
\]%
\[
\varphi \left( \frac{\partial }{\partial z}\right) =\frac{\partial }{%
\partial z}\ ,\qquad \varphi \left( \frac{%
\partial }{\partial t}\right) =\frac{\partial }{\partial t}\ ,\qquad
\varphi \left( \frac{\partial }{\partial s}\right) =0\ ,
\]
\[
g=(dx)^{2}+(dy)^{2}+(dz)^{2}+(dt)^{2}-\eta\otimes \eta.
\]
Then $(\varphi ,\xi ,\eta, g )$\ is a Lorentzian almost paracontact
structure in $ \overline{M}$.

Let $M$ be a hypersurface of $ \overline{M}$ which is defined by
$s=x$ with the imbedding map $i:M\rightarrow \overline{M} $ . Then
the set
\[
\{u_{1}=(1,0,0,0,1),u_{2}=(0,1,0,0,0),u_{3}=(0,0,1,0,0),u_{4}=(0,0,0,1,0)\}
\]
is a local basis for the tangent hyperplane of $M$  and
$N=(1,0,0,0,-1)$ is the normal vector field of the hypersurface.
Since $\xi_{i(p)}=\frac{1}{2}(u_{1}-N)_{i(p)}$, it can be easily
seen that the characteristic vector field $\xi_{i(p)}$, $p\in M$,
does not belong to the tangent hyperplane of $M$. Moreover, since
$\varphi u_{1}=u_{1},\,\varphi u_{2}=u_{2},\,\varphi
u_{3}=u_{3},\,\varphi u_{4}=u_{4}$, then  $M$ is an invariant
hypersurface of $\overline M$ with the characteristic vector field
$\xi_{i(p)}$, $p\in M$, which does not belong to the tangent
hyperplane of the hypersurface.
\end{example}

\begin{example}
Let $ \overline{M}$ be the 3-dimensional real number space with a
coordinate system $(x,y,z)$. If we define
\[
\eta =dz\ ,\qquad \xi =-\,\frac{\partial }{\partial z}\ ,
\]%
\[
\varphi \left( \frac{\partial }{\partial x}\right) =-\,\frac{\partial }{%
\partial x}\ ,\qquad \varphi \left( \frac{%
\partial }{\partial y}\right) =-\,\frac{\partial }{\partial y}\ ,\qquad
\varphi \left( \frac{\partial }{\partial s}\right) =0,
\]
\[
g=(dx)^{2}+(dy)^{2}-\eta\otimes \eta.
\]
on $ \overline{M}$, then $(\varphi ,\xi ,\eta, g )$\ is a Lorentzian
almost paracontact structure in $ \overline{M}$.

Assume that $M$ be a surface of $ \overline{M}$ given by $x=\arcsin
y$ with the imbedding map $i:M\rightarrow \overline{M} $ . Then
\[
\{u_{1}=(1,\sqrt{1-y^{2}},0),u_{2}=(0,0,1)\}
\]
forms a local basis for the tangent plane of $M$  and
$N=(\sqrt{1-y^{2}},-1,0)$ is the normal vector field of the surface.
For any tangent vector field $X$ of the surface we have
\begin{eqnarray}
\varphi i_{*}X=-f_{1}u_{1}\label{ex-2}
\end{eqnarray}
where $X \equiv i_{*}X=f_{1}u_{1}+f_{2}u_{2}$ for some smooth
functions $f_{1},f_{2}$ on $M$. From (\ref{ex-2}) we obtain that $M$
is an invariant surface of $ \overline{M}$ with the characteristic
vector field $\xi_{i(p)}$, $p\in M$, belonging to the tangent plane
of the surface.
\end{example}

\begin{example}
Let $\overline{M}=R^{3}$ be the 3-dimensional real number space with
a coordinate system $(x,y,z)$. We define
\begin{eqnarray}
\eta &=&dz,\,\,\,\,\, \xi =-\frac{\partial }{\partial z}\,,\nonumber\\
\varphi (\frac{\partial }{\partial x}) &=&\frac{\partial }{\partial x}\,,%
\,\,\,\,\varphi (\frac{\partial }{\partial y})=-\frac{\partial }{\partial y}\,,%
\,\,\,\,\,\varphi (\frac{\partial }{\partial z})=0\,,\label{ex-3}\\
 g&=&e^{-2z}(dx)^{2}+e^{2z}(dy)^{2}-(dz)^{2}.\nonumber
\end{eqnarray}
Then $(\varphi,\xi,\eta,g)$ is a Lorentzian para-Sasakian structure
on $\overline{M}$.

Let $M_{1}$ be a surface of $\overline{M}$ with the imbedding map
$i:M_{1}\rightarrow \overline{M} $ which is given by
\begin{eqnarray*}
z=x+y.
\end{eqnarray*}
Then ${u_{1}=(1,0,1),u_{2}=(0,1,1)}$ is a local basis for the
tangent plane of the surface. The vector field
\begin{eqnarray*}
N=(e^{2(x+y)},e^{2(x+y)},1)
\end{eqnarray*}
is a normal vector field of $M_{1}$. Since
\begin{eqnarray*}
\xi=-\frac{1}{e^{2(x+y)}+e^{-2(x+y)}-1}((e^{2(x+y)})u_{1}+(e^{-2(x+y)})u_{2}-N)
\end{eqnarray*}
then for each $p\in M_{1}$ the characteristic vector field
$\xi_{i(p)}$ does not belong to the tangent plane of the surface. A
tangent vector field of the surface can be written by $X \equiv
i_{*}X=f_{1}u_{1}+f_{2}u_{2}$ for some smooth functions
$f_{1},\,f_{2}$ on $M$. By using (\ref{ex-3}) we have
\begin{eqnarray}
\varphi i_{*}X=f_{1}u_{1}-f_{2}u_{2}+(f_{1}-f_{2})\xi.\label{ex-4}
\end{eqnarray}
From (3.1) and (\ref{ex-4}) we get
$$
i_{*}JX=f_{1}u_{1}-f_{2}u_{2}
$$
and
$$
\alpha(X)=f_{1}-f_{2}
$$ where $J$ acts an almost product structure
on $M_{1}$. Consequently, $M_{1}$ is a noninvariant surface of the
Lorentzian para-Sasakian manifold $\overline{M}$ with $\xi$ nowhere
tangent to $M_{1}$.

Let $M_{2}$ be another surface of $\overline{M}$ which is given by
\begin{eqnarray*}
x=\arctan y.
\end{eqnarray*}
Then ${v_{1}=(\frac{1}{1+y^{2}},1,0),v_{2}=(0,0,1)}$ forms a local
orthogonal basis for the tangent plane of the surface and
\begin{eqnarray*}
N=(e^{2z},-\frac{1}{1+y^{2}}e^{-2z},0)
\end{eqnarray*}
is a normal vector field of $M_{2}$. It is obvious that the
characteristic vector field of the manifold belongs to the tangent
plane of $M_{2}$. For any tangent vector field $i_{*}Y\equiv Y$ of
the surface where $i:M_{2}\rightarrow \overline{M} $ is an imbedding
map into the Lorentzian para-Sasakian manifold $\overline M$ we can
write $i_{*}Y=\gamma_{1}v_{1}+\gamma_{2}v_{2}$ for some smooth
functions $\gamma_{1},\,\gamma_{2}$ on $M_{2}$. By using
(\ref{ex-3}) we have
\begin{eqnarray*}
\varphi
i_{*}Y=-\gamma_{1}(v_{1}-\frac{2(1+y^{2})}{(1+y^{2})^{2}e^{2z}-e^{-2z}}N).
\end{eqnarray*}
which shows that $M_{2}$ is a noninvariant surface of the Lorentzian
para-Sasakian manifold $\overline{M}$ with $\xi$ tangent to the
surface.
\end{example}

\bigskip \noindent Correspondence Address: \medskip

\noindent Selcen Y\"{u}ksel Perkta\c{s}

\noindent Department of Mathematics, Faculty of Arts and Sciences, \.{I}n%
\"{o}n\"{u} University

\noindent 44280 Malatya, Turkey

\noindent Email: selcenyuksel@@inonu.edu.tr

\medskip

\noindent Erol K\i l\i \c{c}

\noindent Department of Mathematics, Faculty of Arts and Sciences, \.{I}n%
\"{o}n\"{u} University

\noindent 44280 Malatya, Turkey

\noindent Email: ekilic@@inonu.edu.tr

\medskip

\noindent Sadik Kele\c{s}

\noindent Department of Mathematics, Faculty of Arts and Sciences, \.{I}n%
\"{o}n\"{u} University

\noindent 44280 Malatya, Turkey

\noindent Email: keles@@inonu.edu.
\end{document}